\documentclass[10pt]{amsart}
\usepackage{amsmath, amsthm, amssymb}
\usepackage{amsfonts}
\usepackage{mathrsfs}
\usepackage[all]{xy}
\usepackage{wasysym}
\usepackage{bbm}
\usepackage{textcomp}
\usepackage{verbatim}
\usepackage{charter}
\usepackage{fancyhdr}
\usepackage{epsfig}
\usepackage{pb-diagram,pb-xy} 
\usepackage{pstricks}             
\usepackage{pst-all}               
\usepackage{blkarray}

\usepackage{hyperref}
\hypersetup{
    colorlinks,
    citecolor=blue,
    filecolor=blue,
    linkcolor=blue,
    urlcolor=blue
}

\newcommand{\C}{\mathcal{C}}
\newcommand{\N}{\mathcal{N}}
\newcommand{\U}{\mathcal{U}}

\newcommand{\R}{\mathbb{R}}
\newcommand{\im}{\textnormal{im}}

\newcommand{\dhom}{d_{\textnormal{hom}}}

\theoremstyle{plain}% default
\newtheorem{theorem}{Theorem}[section]

\theoremstyle{definition}
\newtheorem{definition}[theorem]{Definition}
\newtheorem{remark}[theorem]{Remark}
\newtheorem{example}[theorem]{Example}
\newtheorem{lemma}[theorem]{Lemma}

\begin{document}

\title{The Leray Dimension of a Convex Code}

\author{Carina Curto} 
\address{Department of Mathematics, The Pennsylvania State University, State College, PA, 16802, USA}
\email{ccurto@psu.edu}
\author{Ram{\'o}n Vera}
\address{Department of Mathematics,  The Pennsylvania State University, State College, PA, 16802, USA}
\email{rvera.math@gmail.com, rxv15@psu.edu}

\begin{abstract}
Convex codes were recently introduced as models for neural codes in the brain.  Any convex code $\C$ has an associated minimal embedding dimension $d(\C)$, which is the minimal Euclidean space dimension such that the code can be realized by a collection of convex open sets.  In this work we import tools from combinatorial commutative algebra in order to obtain better bounds on $d(\C)$ from an associated simplicial complex $\Delta(\C)$.  In particular, we make a connection to minimal free resolutions of Stanley-Reisner ideals, and observe that they contain topological information that provides stronger bounds on $d(\C)$.  This motivates us to define the {\it Leray dimension} $d_L(\C),$ and show that it can be obtained from the Betti numbers of such a minimal free resolution.  We compare $d_L(\C)$ to two previously studied dimension bounds, obtained from Helly's theorem and the simplicial homology of $\Delta(\C)$.  Finally, we show explicitly how $d_L(\C)$ can be computed algebraically, and illustrate this with examples. 
\end{abstract}

\maketitle
%\tableofcontents

\section{Introduction}
Given a collection $\U = \{U_1,\ldots,U_n\}$ of open sets in a topological space $X$, one can define a binary code,
$$\C(\U) :=\left\lbrace \sigma \subseteq [n] \mid U_\sigma \setminus \bigcup_{j \in [n] \setminus \sigma} U_j \neq \emptyset \right\rbrace,$$
where $[n] = \{1,\ldots,n\}$ and $U_{\sigma} = \cap_{i\in \sigma} U_i$.\footnote{By convention, $U_\emptyset = X$, so that $\emptyset \in \C(\U)$ precisely when the sets in $\U$ fail to cover $X$.}
Conversely, given any binary code $\C \subseteq 2^{[n]}$, there exists an open cover $\U$ such that $\C = \C(\U)$.
If all the $U_i$s can be chosen to be convex subsets of $\mathbb{R}^d$, we say that $\U$ gives a {\it convex realization} of $\C$ in dimension $d$.  We say that $\C$ is a {\it convex code} if it has a convex realization.  The smallest $d$ for which this is possible is called the {\it minimal embedding dimension} of $\C$, and is denoted $d(\C)$. 

Convex codes arise in the context of {\it neural coding}, which is the study of how populations of neurons encode information in the brain.  They have, in some sense, been studied for decades in the neuroscience literature, in systems such as primary visual cortex and hippocampus where neurons often display unimodal (and thus convex) receptive fields.  They were first defined mathematically in \cite{neuro-coding, neural-ring}.  A summary of recent results about convex codes, and other topological aspects of neural coding, can be found in \cite{CEB} and references therein.

For a given code $\C$, the minimal embedding dimension $d(\C)$ can be bounded from below by considering $\Delta(\C)$, the smallest simplicial complex that contains $\C$ (see Section~\ref{sec:comparison}).
Two immediate bounds arise from topological considerations. Here we fix a field $\mathbf{k}$ to compute homology groups. The first bound is the {\em homological dimension},
\begin{equation*}
\dhom (\C) :=  \max \lbrace k \mid H_k(\Delta(\C), \mathbf{k}) \not\simeq 0  \rbrace + 1,
\end{equation*}
obtained simply by requiring that the simplicial homology of $\Delta(\C)$ be compatible with the embedding dimension.  The second bound is given by the 
{\em Helly dimension},  %
\begin{equation*}
d_H(\C) :=  \max \lbrace k \mid \Delta(\C) \textnormal{ has a $k$-dimensional induced simplicial hole} \rbrace,
\end{equation*}
and is obtained by looking at topological obstructions arising from hollow simplices inside $\Delta(\C)$ \cite{Cetal15}.

In this work, we introduce a third bound, the {\em Leray dimension} $d_L(\C)$.  This was motivated by the study of free resolutions of Stanley-Reisner ideals, which are algebraic objects naturally related to the simplicial complexes $\Delta(\C)$. Recently, such free resolutions have also been used to obtain results about convexity of neural codes \cite{Cetal15}, though our focus here is on dimension.  The Betti numbers of a minimal free resolution reveal topological information of the simplicial complex and its subcomplexes. The Leray dimension takes into account all nonzero Betti numbers of such a resolution, while the homological and Helly dimensions can be viewed as stemming from a subset of the Betti numbers. Topologically, the Leray dimension can be defined via the homology groups of all induced subcomplexes $\Delta|_\sigma= \lbrace \tau \in \Delta \mid \tau \subseteq \sigma \rbrace$ inside $\Delta = \Delta(\C)$:
\begin{equation*}
d_L(\C) := \max \lbrace k\in \mathbb{Z} \, \mid \, H_k(\Delta|_{\sigma}, \mathbf{k}) \not\simeq 0 \, \textnormal{ for some } \sigma\subseteq [n]  \,  \rbrace  + 1 .
\end{equation*}
Clearly, $d_L(\C)$ provides a lower bound on the minimal embedding dimension $d(\C)$,
$$d(\C) \geq d_L(\C),$$
just as $d_\mathrm{hom}(\C)$ and $d_H(\C)$ are lower bounds.  Since it takes into account homology of all induced subcomplexes of $\Delta$, $d_L(\C)$ is the best of these three bounds. 

Our main result is that $d_L(\C)$, together with $d_\mathrm{hom}(\C)$ and $d_H(\C)$, can be obtained by calculating a minimal free resolution of the Stanley-Reisner ring $S\slash I_{\Delta(\C)}$.  

\begin{theorem}\label{thm:main}
Let $d_L(\C), d_H(\C)$, and $\dhom(\C)$ be the Leray, Helly, and homological dimensions of a code $\C$, with $\Delta = \Delta(\C)$.  Let $S = \mathbf{k}[x_1, \dots, x_n]$ be a polynomial ring in $n$ variables over the field $\mathbf{k}$. Consider a minimal free resolution of the Stanley-Reisner ring $S/I_{\Delta}$, and denote by $\beta_{i,\sigma}(S/I_\Delta)$ the Betti number at step $i$ with grading $\sigma$ of the free resolution.  Define 
\begin{equation}
R_{i, \sigma} := \begin{cases} 
      |\sigma| - i  & \quad \textnormal{if} \quad  \beta_{i,\sigma} (S\slash I_{\Delta})>0 
      \\
      0 &  \quad \textnormal{otherwise}.
   \end{cases}
\end{equation}
Then the dimensions are  
\begin{enumerate}
\item  $ d_L(\C) = \max_{i, \sigma} \lbrace R_{i,\sigma} \rbrace$

\item $d_H(\C) = \max_{\sigma}  \lbrace R_{1, \sigma} \rbrace $

\item $\dhom(\C) = \max_{i} \lbrace R_{i, 1\cdots 1} \rbrace$ 
\end{enumerate}
\noindent In particular, we have that  
$$ d_L(\C) \geq d_H(\C) \, \textnormal{ and } \, d_L(\C) \geq \dhom (\C).$$  
\noindent Furthermore, for each $q\in \mathbb{N}$ there is a $\C$ such that 
$$d_L(\C) - d_H(\C) \geq q \quad \text{ and } \quad d_L(\C) - \dhom(\C) \geq q.$$
\end{theorem}

The first part of the proof follows from an application of Hochster's formula (see Section \ref{sec:defs-and-proofs}), while the last statement is a consequence of an explicit family of examples given in lemma \ref{lemma:cone}.  It follows from the theorem that the Leray dimension is a straightforward bound to compute, because the Betti numbers can be automatically calculated using current computer algebra programs such as Macaulay2 \cite{Mac}.  %. 

The outline of our paper is as follows. In section  \ref{sec:comparison} we compare the homological, Helly, and Leray dimensions.  In section \ref{sec:defs-and-proofs} we present the definitions and objects of commutative algebra that are relevant to the computability of the Leray dimension, as well as to the proofs of our main results.  Finally, we conclude by computing explicit examples in section \ref{sec:examples}. 

\subsection*{Acknowledgments}
We would like to thank Alexander Kunin for his detailed and helpful comments on an earlier draft of this work. We gratefully acknowledge the support of the Statistical and Applied Mathematical Sciences Institute, under grant NSF DMS-1127914.  This work was partially supported by NSF DMS-1225666/1537228 and NSF DMS-1516881 (to CC).

\section{Comparison of dimension bounds}\label{sec:comparison}
In this section we compare the three dimension bounds with some examples. We show that for cones of cross-polytopes we have that $d_L(\C) > \dhom(\C) =1$ and $d_L(\C) > d_H(\C) = 1$.  We also show that neither $d_H(\C)$ nor $\dhom(\C)$ is better than the other.  We can have $d_H(\C) > \dhom(\C)$ as depicted in example \ref{ex:L26}. We can also find that $\dhom(\C) > d_H(\C)$ for cross-polytopes (see lemma \ref{lemma:cross-polytope}).

We begin by reviewing some notions related to the dimension bounds $\dhom(\C), d_H(\C),$ and $d_L(\C)$.  As noted earlier, these bounds depend only on the simplicial complex $\Delta(\C)$ associated to $\C$.
Recall that an {\it (abstract) simplicial complex} $\Delta \subset 2^{[n]}$ is a set of subsets of $[n]$ such that if $\sigma \in \Delta$ and $\tau \subset \sigma$, then $\tau \in \Delta$.
The elements of $\Delta$ are called {\em faces}, and the dimension of a face $\sigma \in \Delta$ is $|\sigma| - 1$.  
To any code $\mathcal{C}$, we can associate the simplicial complex
\begin{equation}
\Delta(\mathcal{C}):= \left\lbrace \sigma \subseteq [n] \mid \sigma \subseteq c \text{ for some } c \in \mathcal{C} \right\rbrace,
\end{equation}
which is the smallest abstract simplicial complex on $[n]$ containing all elements of $\C$.  

Let $\mathcal{U} = \lbrace U_1, \dots, U_n \rbrace$ be a collection of open sets in a topological space $X$.  The {\em nerve} of $\mathcal{U}$ is the simplicial complex 
$$\N(\mathcal{U}):= \lbrace \sigma \subset [n] \, \mid \, U_{\sigma} \not= \varnothing \rbrace.$$  
A simplicial complex is said to be {\it $d$-representable} if it can be realized as the nerve of a collection of {\it convex} open sets in $\mathbb{R}^d$ \cite{Cetal15, TancerSurvey}.
Note that if $\mathcal{U}$ is a convex realization of a code, $\C = \C(\U)$, then $\N(\U) = \Delta(\C)$.  Thus, if $d = d(\C)$ is the minimal embedding dimension of $\C$, then $\Delta(\C)$ is automatically $d$-representable, via the same collection of convex open sets $\U$.  Bounds on the $d$-representability of $\Delta(\C)$ thus give us immediate lower bounds on the minimal embedding dimension $d(\C)$.  

The homological, Helly, and Leray dimensions are all examples of lower bounds on $d(\C)$ stemming from topological obstructions to $d$-representability of $\Delta(\C)$.  The {\em homological dimension} $d_{\textnormal{hom}}(\C)$ is the largest non-trivial simplicial homology group of  $H_{*} (\Delta(\C), \mathbf{k})$. This dimension captures in particular the topology of $\Delta(\C)$, but it is not sensitive to topological obstructions for the embedding of $\C$  that can come from holes of {\it restrictions} of $\Delta(\C)$. On the other hand, Helly dimension captures simplicial holes of induced subcomplexes, such as the empty triangle $124$ in Figure~\ref{figure:L26}. This connects to the notion of a nerve of a cover and $d$-representability of simplicial complexes, which we introduce next.

As a consequence of Helly's theorem, if a $d$--representable simplicial complex $\Delta$ contains all possible $d$--dimensional faces, then it is the full simplex.   

\begin{theorem}[Helly's theorem \cite{helly}]
Let $\U=\{U_1, \ldots, U_n\}$ be a collection of convex open sets in $\R^d$.  If for every $d+1$ sets in $\U$, the intersection is non-empty, then the full intersection $\bigcap_{i=1}^n U_i \not= \varnothing$.
\end{theorem}

A $d$-representable simplicial complex $\Delta$ does not contain an induced $k$--dimensional simplicial hole for $k \geq d$. A simplicial complex is said to contain an induced {\em $k$--dimensional simplicial hole} if it contains $k+1$ vertices such that the induced subcomplex $\Delta_{k+1}$ is isomorphic to a hollow simplex.  By a {\em hollow simplex} we mean a simplicial complex that contains all subsets except the top-dimensional face. The {\em Helly dimension} of $\Delta(\C)$ is the dimension of the largest induced simplicial hole of $\Delta(\C)$.% 

\subsubsection*{{\bf Clique complexes}}
Unfortunately, the Helly dimension can provide a fairly poor bound on $d(\C)$.  An extreme case is when $\Delta(\C)$ is a clique complex.
Recall that a {\it clique} in a graph $G$ is an all-to-all connected subset of vertices in $G$. Note that if $\sigma$ is a clique of $G$, then all subsets of $\sigma$ are also cliques.  The set of all cliques of $G$ is thus naturally a simplicial complex, called the {\it clique complex} $X(G)$:
$$ X(G) = \lbrace \sigma \subseteq [n]  \mid \sigma \text{ is a clique of }  G \rbrace. $$
Because, by definition, a clique complex has no induced simplicial holes other than missing edges, the Helly dimension can be at most 1.  This is, in fact, a defining property for clique complexes.

\begin{lemma}\label{lemma:cliqueHelly}
$\Delta$ is a clique complex if and only if $d_H(\Delta) \leq 1$.
\end{lemma}

\subsubsection*{{\bf Leray dimension}}
We now turn to the {\it Leray dimension}.  Recall that the Leray dimension of a code $\C$ with simplicial complex $\Delta = \Delta(\C)$ is 
\begin{equation}
d_L(\C) := \max \lbrace k\in \mathbb{Z} \, \mid \, H_k(\Delta|_{\sigma}, \mathbf{k}) \not\simeq 0 \, \textnormal{ for some } \sigma\subseteq [n]  \,  \rbrace  + 1 
\end{equation}
This is closely related to well-known concept of Leray number.
\begin{definition}[\cite{TancerSurvey}]
Let $\Delta$ be a simplicial complex. $\Delta$ is {\em d-Leray} if for every $k\geq d$ we have that the reduced homology groups $\tilde{H}_k (\Delta|_{\sigma}, \mathbf{k}) \simeq 0$ for every induced subcomplex $\Delta|_{\sigma}$.  The {\em Leray number} is the smallest possible $d$ such that $\Delta$ is $d$-Leray. 
\end{definition}
The Leray dimension of a code $\C$ is thus equal to the Leray number of $\Delta(\C)$ plus one:
\begin{equation}
d_L(\C) = \textnormal{ Leray number of } \Delta(\C) + 1.
\end{equation}
%. 

\begin{remark}
Although it is known that the minimal embedding dimension $d(\C)$ depends on details of the code beyond $\Delta(\C)$ \cite{neural-ring}, all dimensions we consider in this work depend only on $\Delta(\C)$. We will thus write $d_L(\Delta), d_H(\Delta), \dhom(\Delta)$ to denote the corresponding dimensions for any code $\C$ with simplicial complex $\Delta = \Delta(\C)$. 
\end{remark}

\begin{remark}
The dimension bounds considered in this work are defined with respect to homology groups that depend on a field $\mathbf{k}$. By abuse of notation we will suppress the field $\mathbf{k}$ in our notation, and simply write $d_L(\Delta), d_H(\Delta),$ and $\dhom(\Delta)$. 
\end{remark}

The next example illustrates a case where $d_L(\Delta) = d_H(\Delta) > \dhom(\Delta)$. 
\begin{figure}[h]
\begin{center}
\includegraphics[width=5in]{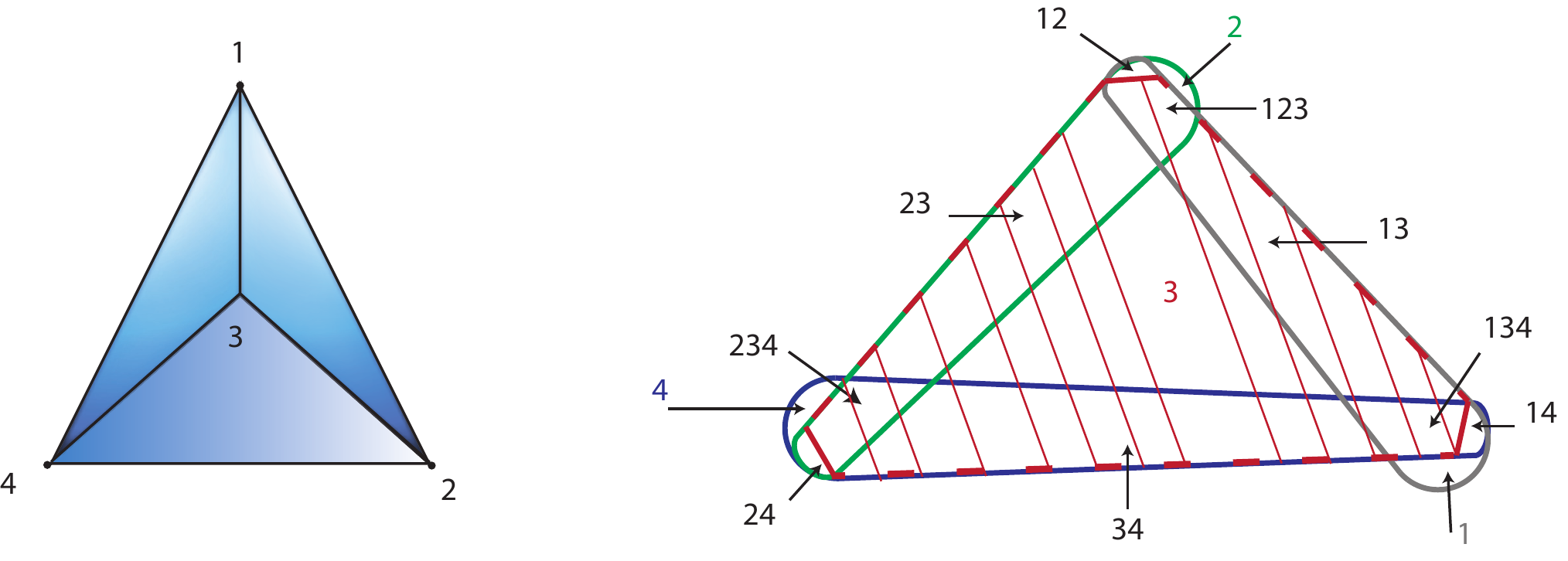}
\end{center}
\caption{(Left) The simplicial complex $\Delta(\C)$ for the code in example \ref{ex:L26}. Note that the triangle $124$ is missing. (Right) A convex realization of $\C$ in $\mathbb{R}^2$.}
\label{figure:L26}
\end{figure}

\begin{example}\label{ex:L26}
Let $\C$ be the code on $n=4$ neurons, with the following codewords:
$$\lbrace  0000, 1000, 0100, 0010, 0001,  1100, 1010, 1001, 0110, 0101, 0011, 1110,  1011, 0111 \rbrace.$$
The simplicial complex $\Delta = \Delta(\C)$ is shown in Figure~\ref{figure:L26} (left). Notice that this is not a clique complex, as the face $124$ is missing. The restriction $\Delta|_{\{1,2,4\}}$ reveals an induced simplicial hole of dimension 2, thus $d_H(\Delta) = 2$. Looking at the homology of all induced subcomplexes we see that the Leray dimension is also $d_L(\Delta) = 2$.  On the other hand, the whole simplicial complex $\Delta$ is contractible, so $\dhom(\Delta) = 1$.    Since the minimal embedding dimension $d(\C) \geq d_L(\Delta) = 2$, the existence of a convex realization of $\C$ in dimension $2$ shows that $d(\C) = 2$ (see Figure~\ref{figure:L26}, right) .

\end{example}

The following family of polytopes has $d_L(\Delta) = \dhom(\Delta) > d_H(\Delta)$. 

\subsubsection*{{\bf Cross-polytopes}}\label{sec:cross-polytopes}
\begin{figure}[h]
\begin{center}
\resizebox{.99\textwidth}{!}{
\includegraphics{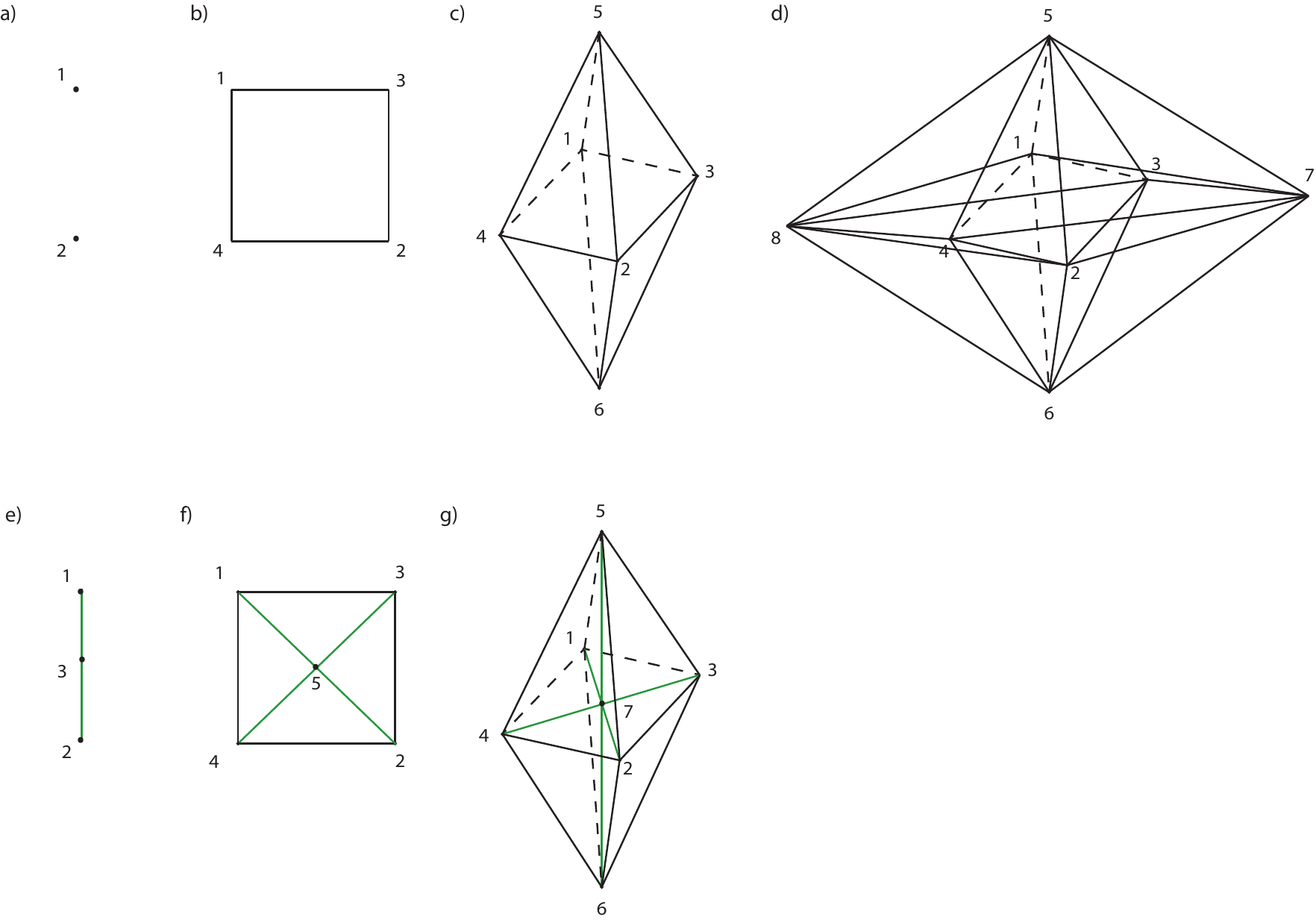} }
\end{center}
\caption{The first four cross-polytopes are: (a)  two vertices $\Gamma_0$, (b) a square $\Gamma_1$, (c) an octahedron $\Gamma_2$,  and (d) an orthoplex $\Gamma_3$ (all faces corresponding to cliques are filled in). The first three cones of cross-polytopes are:  (e) the cone over two vertices $\hat{\Gamma}_0$,  (f) the cone over the square $\hat{\Gamma}_1$, and (g) the cone over the octahedron $\hat{\Gamma}_2$. }
\label{fig:c-polytopes}
\end{figure}

Consider a graph $G_0$ with two vertices and no edges as in Figure \ref{fig:c-polytopes}a. Add two new points and connect each of the two vertices of the previous graph with the new vertices. This gives rise to a new graph, $G_1$, that looks like a square (Figure \ref{fig:c-polytopes}b). To construct the next graph $G_2$, add two new vertices to $G_1$ and connect them to the previous four vertices.  We can proceed inductively to obtain the sequence of graphs $G_0, G_1, G_2, \ldots$. By taking clique complexes of these graphs, we obtain a family of cross-polytopes
$$\Gamma_i = X(G_i).$$
Starting with $G_2$, this process fills in higher-dimensional faces to obtain an octahedron $\Gamma_2$ (Figure \ref{fig:c-polytopes}c), an orthoplex $\Gamma_3$ (Figure \ref{fig:c-polytopes}d), and so on.  

The Helly dimension of any cross-polytope $\Gamma_i$ is $d_H(\Gamma_i) = 1$, since they are all clique complexes (see lemma \ref{lemma:cliqueHelly}).  Yet, the homological dimension is $\dhom(\Gamma_i) = i + 1$, since it detects the largest non-trivial homology group.  We record these observations in the following lemma.

\begin{lemma}\label{lemma:cross-polytope}
Let $\Gamma_i$ be the $i$-th cross-polytope. Then $d_H(\Gamma_i) = 1$ and $\dhom(\Gamma_i) = i +1$. 
\end{lemma}

\subsubsection*{{\bf Cones of cross-polytopes}}\label{sec:cones}
We now consider cones over the same family of cross-polytopes, $\lbrace \Gamma_i \rbrace$.  For each cross-polytope,  we add a vertex in the middle and join all edges and faces to form a cone.  This produces a new family of simplicial complexes, $\lbrace \hat{\Gamma}_i\rbrace$, which we call the cones of the cross-polytopes. Figure \ref{fig:c-polytopes}e-g depicts the first three such cones: $\hat{\Gamma}_0$, $\hat{\Gamma}_1$, and $\hat{\Gamma}_2$.  Note that $ \hat{\Gamma}_i$ is contractible for any $i$, so $\dhom(\hat{\Gamma}_i) = 1$.
On the other hand, the Leray dimension is $d_L(\hat{\Gamma_i}) = i + 1$ because of the homology of the subcomplex obtained by removing the added cone point.  As before, the Helly dimension is $d_H (\hat{\Gamma}_i) = 1$.  We collect these observations in the following lemma. 

\begin{lemma}\label{lemma:cone}
Let $\Gamma_i$ be the $i$-th cross-polytope and $\hat{\Gamma}_i$ the cone over $\Gamma_i$. Then $d_H(\hat{\Gamma}_i) =1$, $\dhom(\hat{\Gamma}_i) = 1$, but $d_L(\hat{\Gamma}_i) = i +1 $. 
\end{lemma}

%
%
%%%%%%%%%%%%%%%%%%%%%%%%%%%%%%%%%%%%%%%%%
%
%  SECTION 2  DEFINITIONS AND PROOFS
%
%%%%%%%%%%%%%%%%%%%%%%%%%%%%%%%%%%%%%%%%%
%
\section{Computing Leray dimension algebraically}\label{sec:defs-and-proofs}
In this section we recall the topological and combinatorial concepts that are relevant to prove the theorem of this work. In particular, we go through the definitions of Stanley-Reisner ideal, free resolutions, Betti numbers and Hochster's formula. %

\subsection{Stanley-Reisner ideal and free resolutions}
We follow \cite{MSbook} for the next definitions. Recall that a monomial $x_1^{a_1}\cdots x_n^{a_n}$ is {\em squarefree} if every exponent $a_i$ is 0 or 1. An ideal is {\em squarefree} if it is generated by squarefree monomials. There is a bijective correspondence between simplicial complexes on $n$ vertices and squarefree monomial ideals on $x_1,\ldots,x_n$ \cite{MSbook}. In what follows, we denote squarefree monomials as $x_{\sigma} = \prod_{i\in \sigma} x_i$. 

\begin{definition}
Let $\Delta$ be a simplicial complex. The {\em Stanley-Reisner ideal} of $\Delta$ is the squarefree monomial ideal
$$I_{\Delta} = \langle x_{\sigma} \mid \sigma \notin \Delta \rangle. $$
Its generators are monomials corresponding to nonfaces $\sigma$ of $\Delta$.  The {\em Stanley-Reisner ring} of $\Delta$ is the quotient ring  $S\slash I_{\Delta},$ where $S = \mathbf{k}[x_1, \dots, x_n]$ denotes the polynomial ring over a field $\mathbf{k}$. 
\end{definition}  
 
We are interested in resolutions of Stanley-Reisner ideals, since these are related to simplicial complexes.  The direct sum module $F \simeq S\oplus \dots \oplus S$ of $r$ copies of $S$ is called the {\em free $S$-module of finite rank $r$}. 
Here we consider $S$ to be $\mathbb{N}^n$--graded, meaning that $S\simeq S(-{\bf a_1}) \oplus \dots \oplus S(-{\bf a_r})$ for some vectors ${\bf a_1, \dots , a_r} \in \mathbb{N}^n$.  A sequence of maps of free $S$-modules 
$$ 0 \longleftarrow  F_0 \stackrel{\phi_1}\longleftarrow  F_1 \longleftarrow \ldots \longleftarrow F_{k-1} \stackrel{\phi_k}\longleftarrow F_{k} \longleftarrow 0 $$
is a {\em chain complex} if $\phi^2 = 0$.  A chain complex is {\em exact} in homological degree $k$ if $\ker(\phi_k) = \im(\phi_{k+1})$. %. 
\begin{definition}
Let $S = \mathbf{k}[x_1, \dots, x_n]$ be a polynomial ring in $n$ variables over a field $\mathbf{k}$. A {\em free resolution} of an ideal $I$ of $S$ is an exact sequence of free modules
\begin{equation}\label{def:resolution}
0 \longleftarrow S\slash I \longleftarrow F_0 \stackrel{\phi_1}\longleftarrow  F_1 \longleftarrow \ldots \longleftarrow F_{i-1}  \stackrel{\phi_i}\longleftarrow F_{i} \longleftarrow \ldots  \stackrel{\phi_m}\longleftarrow F_\ell \longleftarrow 0 
\end{equation}
The {\em length of the resolution} is the greatest homological degree $\ell$ of a nonzero module in the resolution.\end{definition}
\noindent Every finitely-generated module over a polynomial ring has a free resolution of finite length. In our context, $S\slash I$ is $\mathbb{N}^n$--graded, thus it has an $\mathbb{N}^n$--graded free resolution.

\begin{remark}\label{rem:notation-F}
Throughout this work we will abuse notation and denote the modules of a free resolution as $F_{i}^{k}$, where $k$ denotes the number copies of $S$ and the subindex $i$ corresponds to the step of the resolution. 
\end{remark}

An explicit way to see the maps between free modules is through the matrix representing them. Denote by $\succeq$ the partial order on $\mathbb{N}^n$ where $\mathbf{a} \succeq \mathbf{b}$ if and only if $a_i \geq b_i$ for all $i \in [n]$.  A {\em monomial matrix} is an array of scalar entries $\lambda_{qp}$ whose columns correspond to the source degrees $\mathbf{a_p}$  and whose rows represent the target degrees $\mathbf{a}_q$.  The scalar entry indicates that the basis vector of $S (- {\bf a}_p )$ should map to an element that has coefficient $\lambda_{qp}$ on the monomial. That is ${\bf x}^{{\bf a}_p -  {\bf a}_q} $ times the basis vector of  $S (- {\bf a}_q )$.  %.  
The requirement $ {\bf a}_p \succeq {\bf a}_q$ guarantees that $\bf {x}^{{\bf a}_p -  {\bf a}_q} $ has nonnegative exponents.   The monomial matrix of a free resolution representing $\phi_i$  has the following form 

$$\bigoplus_q S(-{\bf x^{a_q}})   \xleftarrow{\left[ \begin{array}{cccc}
                &  \dots  &  {\bf a}_ p  & \dots  \\
 \vdots      &          &        &  \\
 {\bf a}_q  &   &   \left[ {  \bf x}^{{\bf a}_p -  {\bf a}_q}   \lambda_{qp} \right]  &    \\
 \vdots    &          &        &
\end{array} \right] }   \bigoplus_p S(-{\bf x^{a_p}} ) $$

\begin{definition}
A monomial matrix is {\em minimal} if $\lambda_{qp} = 0$ when $ {\bf a}_p = {\bf a}_q $. A free resolution is {\em minimal} if it can be written down with minimal monomial matrices. 
\end{definition}

%%%%%%%%%%%%%%%%%%%%%%%%%%
% BETTI NUMBERS AND HOCHSTER'S FORMULA
%%%%%%%%%%%%%%%%%%%%%%%%%%

\subsection{Betti Numbers and Hochster's Formula}
First, we recall Hochster's formula, which relates the Betti numbers of a minimal free resolution of $I_{\Delta}$ to topological invariants of $\Delta$ and its subcomplexes. 
\begin{definition}
Let $\mathcal{F}$ be a minimal free resolution of a finitely $\mathbb{N}^n$--graded module $S/I$.   The {\em Betti number} $\beta_{i,\sigma}(S\slash I)$ is the rank of the module in multidegree $\sigma$ at step $i$ of the resolution $\mathcal{F}$, where $S/I$ is step 0 and the steps increase as we move from left to right.
\end{definition}

\noindent Hochster's formula focuses on minimal free resolutions of Stanley-Reisner ideals $I_\Delta$.  

\begin{theorem}[{Hochster's Formula \cite[Corollary 5.12]{MSbook}}]\label{thm:Hochsters}
Let $S = \mathbf{k}[x_1, \dots, x_n]$ be a polynomial ring, $\Delta$ a simplicial complex, $I_\Delta$ its Stanley-Reisner ideal, and $\beta_{i,\sigma}(S\slash I_{\Delta})$ the Betti numbers of a minimal free resolution of $S\slash I_{\Delta}$.  The nonzero Betti numbers lie only in squarefree degrees $\sigma$, and we have
\begin{equation}
\beta_{i, \sigma} (S\slash I_{\Delta}) = \dim_{\mathbf{k}} \tilde{H}^{|\sigma| - i - 1} (\Delta|_{\sigma}, \mathbf{k})
\end{equation}
\end{theorem}
 
\noindent This allows us to compute the Leray dimension from minimal free resolution, as illustrated in the next example.

\begin{example}
Consider the simplicial complex $\Delta = \Gamma_1$ in Figure \ref{fig:c-polytopes}b. The Stanley-Reisner ideal is given by 
$$I_{\Delta} = \langle x_1 x_2, x_3 x_4 \rangle $$
Computing the minimal free resolution of $S\slash I_{\Delta}$, we obtain  
$$   \begin{tabular}{cccccc} 
$0  \leftarrow S\slash I_{\Delta}    \xleftarrow{\makebox[1cm]{}}$ & $F_1^2    \xleftarrow{\makebox[1cm]{}}$            & $F_2  \xleftarrow{\makebox[1cm]{}}$     	  & 0  &   \\   \\ 
											& $\sigma_1 = 1100$				& $\sigma_3 = 1111$     			  &  				  &     \\                                                                                                                                                                                                                                                                                                                                                                                                                                                           											& $\sigma_2 = 0011$      		        		& 	         					  &       			  &  \\
\end{tabular} $$
The Betti numbers are $\beta_{1, \sigma_1} = 1,  \beta_{1, \sigma_2} = 1,  \beta_{2, \sigma_3} = 1$. Note that $|\sigma_1 |= 2, |\sigma_2 | = 2$, both lying at level $i=1$, and also $|\sigma_3 | = 4$ at level $i=2$. 
Using Hochster's formula, we compute:
$$\beta_{1, \sigma_1} =  \beta_{1, \sigma_2} = 1  \Rightarrow  \dim \tilde{H}^{0}(\Delta|_{\sigma_1}, \mathbf{k}) = 1, $$
$$\beta_{2, \sigma_3} = 1  \Rightarrow  \dim \tilde{H}^{1}(\Delta|_{\sigma_3}, \mathbf{k}) = 1. $$
\end{example}

%%%%%%%%%%%%%%%%%%%%%%%%%%%%%%%%%%%%
% PROOF OF THEOREM 1.4 (MAIN THM)
%%%%%%%%%%%%%%%%%%%%%%%%%%%%%%%%%%%%

\subsection{Proof of Theorem \ref{thm:main}}\label{sec:proof2}
The last statement of the theorem was a consequence of lemma \ref{lemma:cone}. Recall that we consider $\Delta = \Delta(\C)$ and 
\begin{equation}
R_{i, \sigma} := \begin{cases} 
      |\sigma| - i  & \quad \textnormal{if} \quad  \beta_{i,\sigma} (S\slash I_{\Delta})>0  \quad \textnormal{and,}
      \\
      0 &  \quad \textnormal{otherwise}
   \end{cases}
\end{equation}
\vspace{2mm}

\noindent {\bf 1.}  $ d_L(\Delta) =  \max_{i, \sigma}\lbrace R_{i,\sigma} \rbrace$.
\\
Recall that the Leray dimension of a simplicial complex $\Delta$ is defined as 
$$d_L(\Delta) = \max \left\lbrace k\in \mathbb{Z} \, \mid \, H_k(\Delta|_{\sigma}, \mathbf{k}) \not\simeq 0 \, \textnormal{ for some } \sigma\subseteq [n]  \, \right\rbrace  + 1 $$
By the universal coefficient theorem, the dimensions of homology groups $H_{*}(\Delta, \mathbf{k})$ correspond to dimensions of cohomology groups $H^{*}(\Delta, \mathbf{k})$ \cite[Thm. 3.2, p 195]{Hatcher}. Hochster's formula \ref{thm:Hochsters} provides a link between cohomology and Betti numbers. Consequently it makes sense to express $d_L(\Delta)$ in terms of Betti numbers. In particular, Hochster's formula gives us a relation between the Betti numbers of the free resolution of the ring $S$ with the Stanley-Reisner ideal $I_{\Delta}$ and the dimension of the cohomology groups of the simplicial complex $\Delta$ restricted to degrees $\sigma$: 
\begin{equation}\label{eq:Hochster}
\beta_{i, \sigma} (S\slash I_{\Delta}) = \dim_{\mathbf{k}} \tilde{H}^{|\sigma| - i - 1} (\Delta|_{\sigma}, \mathbf{k}).
\end{equation}
To obtain the desired result we maximize over all the cohomology groups of the restrictions $\Delta|_{\sigma}$ using the relation of the previous equation. Notice that such a calculation can be performed with a minimal free resolution. Thus we can obtain 
\begin{equation}
d_L(\Delta) = \max_{\beta_{i,\sigma} (S\slash I_{\Delta})>0} (|\sigma| - i ).
\end{equation}
\vspace{2mm}

\noindent {\bf 2.}  $d_H(\Delta) = \max_{\sigma}  \lbrace R_{1, \sigma} \rbrace $.  %
\\
\noindent Denote by $\min I_{\Delta}$ the set of minimal monomial generators of the ideal $I_{\Delta}$.  First observe that $ d_H = \max_{x_\sigma \in \min I_{\Delta}} |\sigma| - 1$ since the minimal monomials correspond to simplicial holes.  Let $\mathcal{F}$ be a minimal free resolution.  On the first level of the resolution, the minimal monomial generators correspond to the elements of the matrix of the first map $\phi_1$.  That is, the minimal monomials give rise to gradings $\sigma$ at level $F_1$ by assigning a 1 to the $i$-th slot of $\sigma$ for every $x_i$ in the monomial and 0 otherwise.  The source elements of $F_1$ are mapped to $S\slash I_{\Delta}$  with $\phi_1$.  This means that $|\sigma|$ for elements $\sigma$ in $F_1$ is precisely the degree of the monomials $x_{\sigma} \in \min I_{\Delta}$.  Thus, $d_H(\Delta)$  is in turn equal to $\max_{\beta_{1,\sigma} (S\slash I_{\Delta})>0} (|\sigma| - 1 )  = \max_{\sigma}\lbrace R_{1, \sigma} \rbrace$, since the Betti numbers $\beta_{1,\sigma} (S\slash I_{\Delta})$ at level 1 are determined by the length of $x_{\sigma} \in \min I_{\Delta}$. Hence, Helly dimension can be computed at the first level of the resolution. 
\vspace{2mm}

\noindent {\bf 3. }  $\dhom(\Delta) = \max \lbrace R_{i, 1\cdots 1} \rbrace$.  %
\\
Since $\Delta|_{\sigma = 1\cdots 1} = \Delta$, then $\dhom(\Delta)$ can be obtained form the resolution by looking at all Betti numbers $\beta_{i, \sigma}$ with $\sigma=1\cdots 1$,  and
$$\dhom(\Delta) = \max _{\beta_{i, 1\cdots 1} (S\slash I_{\Delta})>0} (|\sigma| - i ).$$
%

%
%%%%%%%%%%%%%%%%%%%%%%%%%%%%%%%%%%%%
%
% SECTION         EXAMPLES
%
%%%%%%%%%%%%%%%%%%%%%%%%%%%%%%%%%%%%
%
\section{Examples with Calculations of the Leray Dimension}\label{sec:examples}
In this section we use Theorem \ref{thm:main} to compute the dimension bounds for a variety of examples. Recall that in Example \ref{ex:L26} we had $d_L = d_H > \dhom$. The next example has $d_L = \dhom > d_H$. 
%
%   OCTAHEDRON
%
\begin{example}{\bf (Octahedron)} \label{ex:octahedron} 
Consider a code $\C$ whose simplicial complex $\Delta = \Delta(\C)$ is the octahedron in Figure \ref{fig:c-polytopes}c with 8 faces, 12 edges, 6 vertices, and it is empty inside.  This is an example of a cross-polytope.  The Stanley-Reisner ideal of this simplicial complex is given by 
$$I_\Delta = \langle x_1x_2, x_3 x_4, x_5 x_6 \rangle$$
\noindent A minimal free resolution of $S/I_{\Delta}$ is:
\begin{eqnarray}\label{ResOctahedron}
 0 \longleftarrow 
 S/I_{\Delta} \xleftarrow{[\begin{array}{ccc}x_1x_2 & x_3 x_4 & x_5 x_6\end{array}]} 
 S(-2)\oplus S(-2) \oplus S(-2)  \xleftarrow{\left[\begin{array}{ccc}x_3 x_4 & x_5 x_6 & 0\\ -x_1x_2 & 0 & -x_5 x_6 \\ 0 & -x_1x_2 & x_3 x_4 \end{array}\right]} &&
\nonumber \\ 
 S(-4)\oplus S(-4) \oplus S(-4) 
 \xleftarrow{\left[\begin{array}{c}x_5 x_6\\-x_3x_4\\x_1x_2\end{array}\right]} S(-6) \longleftarrow  0 && 
\end{eqnarray}
The entries of the monomial matrices express the terms necessary to go from a source basis element to a target one. To see this we can write on top of each column the source monomial and next to each row the target. The matrix entry denotes the monomial required to go from one to another. If the monomials of the source and target do not share any common elements then we write a 0 for the corresponding matrix entry. To see this, we take for instance the second monomial matrix and label the rows and columns.
\[
\begin{blockarray}{cccc}
&  x_1x_2x_3x_4 &     x_1x_2x_5x_6    &    x_3 x_4 x_5 x_6     \\
\begin{block}{c[ccc]}
  x_1 x_2    & x_3 x_4 & x_5 x_6 & 0  \\
  x_3 x_4    & -x_1x_2 & 0 & -x_5 x_6 \\
  x_5 x_6    & 0 & -x_1x_2 & x_3 x_4   \\
\end{block}
\end{blockarray}
 \]
 
We use the condensed notation from remark \ref{rem:notation-F} to express the resolution \eqref{ResOctahedron} in terms of modules $F_i^k$, which are the direct sum of $k$ copies of $S$ at level $i$. We write below each module the list of the degrees of its generators as binary patterns.
 \begin{eqnarray}
 0 \xleftarrow{\makebox[1cm]{}}  S/I_{\Delta} \xleftarrow{\makebox[1cm]{}}F_{1}^{3}  \xleftarrow{\makebox[1.5cm]{}} F_2^{3}  \xleftarrow{\makebox[1.5cm]{}}F_3 \longleftarrow  0 
 \nonumber 
\end{eqnarray}
$   \begin{tabular}{ccccc} 
\hspace{5mm}& \hspace{52mm}$\sigma_1 = 110000$		&\hspace{5mm} $\sigma_4 = 111100$     		&  \hspace{3mm}$\sigma_7 = 111111$   \\                                                                                                                                                                                                                                                                                                                                                                                                                                                           			& \hspace{52mm} $\sigma_2 = 001100$         	&\hspace{5mm}  $\sigma_5 = 110011$	         	&      			 \\ 
			& \hspace{52mm} $\sigma_3 = 000011$  	        	& \hspace{5mm}  $\sigma_6 = 001111$			&       		 \\
\end{tabular} $

\noindent The length of these elements is  
$$ |\sigma_k| = 2 \quad \textnormal{for } k=1,2,3, \qquad  |\sigma_{\ell}| = 4 \quad \textnormal{for } \ell=4,5,6,  \quad \textnormal{and} \quad  |\sigma_7| = 6. $$
Recall that the Leray dimension is $ d_L(\Delta) = \max\lbrace R_{i,\sigma} \rbrace$, where 
\begin{equation*}
R_{i, \sigma} := \begin{cases} 
      |\sigma| - i  & \quad \textnormal{if} \quad  \beta_{i,\sigma} (S\slash I_{\Delta})>0 
      \\
      0 &  \quad \textnormal{otherwise}
   \end{cases}
\end{equation*}
Since the Betti numbers $\beta_{i,\sigma}(S/I_{\Delta})$ correspond to the rank of the module in degree $\sigma$ at step $i$, we can directly read from the previous resolution the values of $R_{i, \sigma_{i}}$.  
$$\begin{tabular}{rclcl}
$\beta_{1, \sigma_k} = 1 $  &  $\Rightarrow$  & $ R_{1, \sigma_k} = |\sigma_k| - 1$ = 1 & \quad  & $\text{for } k=1,2, 3, $ 
\\
$\beta_{2, \sigma_{\ell}} = 1 $  &  $\Rightarrow$  & $ R_{2, \sigma_{\ell}}  = |\sigma_{\ell}| - 2 = 2$ &             &$\text{for } \ell= 4,5,6,$ 
\\
$\beta_{3, \sigma_7} = 1 $  &  $\Rightarrow$  & $ R_{3, \sigma_7} = |\sigma_7| - 3 = 3 $ &       &
\end{tabular}$$
Applying Theorem \ref{thm:main} we obtain 
\begin{itemize}
\item[]  $ d_L(\Delta) = \max_{i, \sigma} \lbrace R_{i,\sigma} \rbrace = 3$,

\item[] $d_H(\Delta) = \max_{\sigma}  \lbrace R_{1, \sigma} \rbrace = 1$,

\item[] $\dhom(\Delta) = \max_{i} \lbrace R_{i, 1\cdots 1} \rbrace = 3$.
\end{itemize}
Notice that Helly and homological dimension can also be seen from the topology of $\Delta$. Since $\Delta$ is a clique complex then it follows from lemma \ref{lemma:cliqueHelly} that $d_H(\Delta) = 1$, and since the octahedron has the homology of the 2-sphere, we have $\dhom(\Delta) = 3$.  
\end{example}

%
% CONE OF A CROSS-POLYTOPE
%
\begin{example}{\bf(Cone of a Cross-polytope)}\label{ex:envelope}
This example shows $d_L >  d_H > \dhom$. Consider a code $\C$ having a simplicial representation as in Figure \ref{fig:c-polytopes}f. Recall that this is a cone of a cross-polytope. 

\begin{figure}[h]\label{fig:envelope}
\begin{center}
\resizebox{.75\textwidth}{!}{
\includegraphics{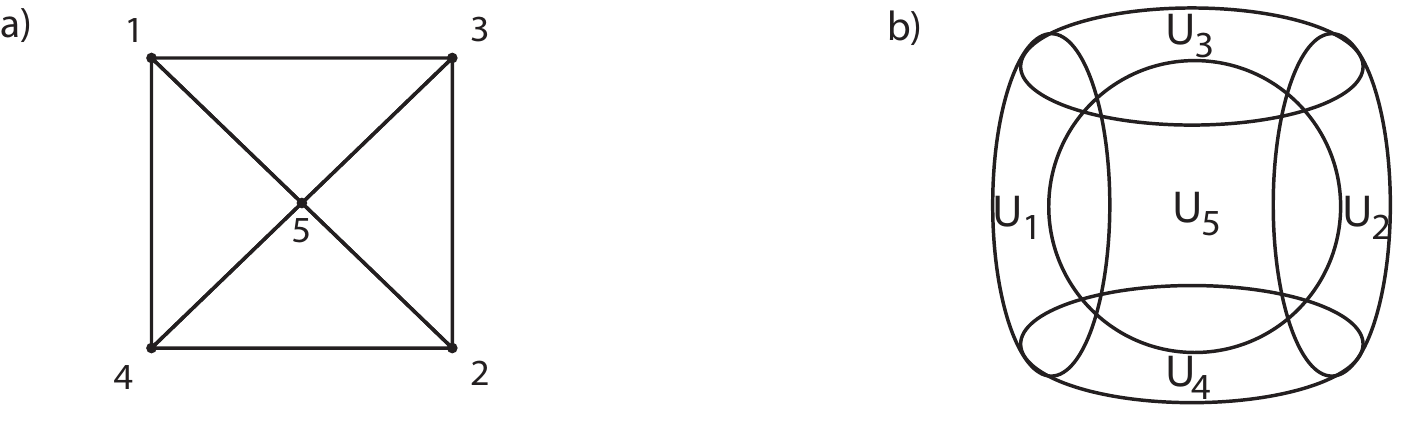} }
\end{center}
\caption{(a) Simplicial complex $\Delta$ of code $\C$ of example \ref{ex:envelope}. The faces are filled in.  (b) A convex realization of $\C$ in $\mathbb{R}^2$.}
\label{fig:envelope}
\end{figure}
 
The Stanley-Reisner ideal of this simplicial complex is   $I_{\Delta} = \langle x_1 x_2, x_3 x_4 \rangle$.  A minimal free resolution of $I_{\Delta}$ is,

\begin{eqnarray}
 0 \xleftarrow{\makebox[1cm]{}}  S/I_{\Delta}   \xleftarrow{\makebox[1.5cm]{}} F_1^{2}  \xleftarrow{\makebox[1.5cm]{}}F_2 \xleftarrow{\makebox[1cm]{}}  0 ,
 \nonumber 
\end{eqnarray}
$   \begin{tabular}{ccccc} 
\hspace{5mm}& \hspace{60mm}$\sigma_1 = 11000$	&\hspace{5mm} $\sigma_3 = 11110$ 		&    \\                                                                                                                                                                                                                                                                                                                                                                                                                                                           			& \hspace{60mm} $\sigma_2 = 00110$     &        	&      			 \\ 

\end{tabular} $
\\

\noindent with non-zero Betti numbers,
$$\beta_{1, \sigma_1} = 1 \quad , \quad \beta_{1, \sigma_2} = 1 \quad , \quad \beta_{2, \sigma_3} = 1,$$
leading to 
$$R_{1, \sigma_{1}} = |\sigma_1| - 1 =  1  \quad , \quad R_{1, \sigma_{2}} = |\sigma_2| - 1 = 1 \quad , \quad R_{2, \sigma_{3}} =  |\sigma_3| - 2  =2$$
It follows from Theorem \ref{thm:main} that
\begin{itemize}
\item[]  $ d_L(\Delta) = \max_{i, \sigma} \lbrace R_{i,\sigma} \rbrace = 2$

\item[] $d_H(\Delta) = \max_{\sigma}  \lbrace R_{1, \sigma} \rbrace = 1$

\item[] $\dhom(\Delta) = \max_{i} \lbrace R_{i, 1\cdots 1} \rbrace = 0.$ 
\end{itemize}

\end{example}

\medskip

%
%   CALCULATING WITH MACAULAY2
%

We now illustrate how the dimensions can be obtained using the computer algebra package Macaulay2 \cite{Mac}. First we give the general framework assuming that the Stanley-Reisner ideal of $\Delta$ is given. 

\subsubsection*{General Framework}
\quad
\\
\noindent {\bf Step 1}: Define the algebraic objects to be worked upon: ground field $\mathtt{kk}$, polynomial ring $\mathtt{R}$, Stanley-Reisner ideal $\mathtt{SRideal}$, and module $\mathtt{M}$
\medskip

\noindent
$\mathtt{kk = ZZ/2;}$
\medskip 

\noindent
$\mathtt{R = kk[ x_1, \dots , x_n, Degrees => \lbrace \lbrace 1, 0, \dots, 0\rbrace, \lbrace 0,1, \dots, 0 \rbrace, \dots \lbrace 0, \dots, 0, 1 \rbrace \rbrace }]; $
\medskip

\noindent
$\mathtt{SRideal = monomialideal( x_{\sigma_{1}}, \dots , x_{\sigma_{\ell}}) }$ 
\medskip

\noindent
$\mathtt{M = R\,\hat{}\, 1 / SRideal }$
\medskip

\noindent {\bf Step 2}: Calculate the minimal free resolution
\medskip

\noindent $\mathtt{Mres = res\, \, M;}$
\medskip 

\noindent {\bf Step 3}: Obtain the Betti numbers coming from the free resolution
\medskip

\noindent $\mathtt{peek\, \, \, betti \, \, \,  Mres}$
\medskip

\noindent It follows from theorem \ref{thm:main} that the Leray dimension of a code $\C$ is $d_L(\Delta) = \max_{\beta_{i,\sigma} (S\slash I_{\Delta})>0} (|\sigma| - i ) $

%\subsubsection*{Example}
\bigskip

\noindent We proceed a concrete example. Consider the simplicial complex of example \ref{ex:envelope}. 
\medskip

\noindent $\mathtt{i1 : kk=ZZ/2;}$
\medskip

\noindent  $\mathtt{i2 : R = kk[x1,x2,x3,x4,x5, Degrees => \lbrace \{1,0,0,0,0\} , \{0,1,0,0,0\} , \{0,0,1,0,0\} , \{0,0 ,0,1,0\} , \{ 0,0,0,0,1\} \rbrace ] } ;$
\medskip

\noindent $\mathtt{o2 = R}$
\medskip

\noindent $\mathtt{ o2 : PolynomialRing}$
\medskip

\noindent $\mathtt{i3 : SRideal = monomialIdeal(x1*x2, x3*x4);}$ 
\medskip

\noindent $\mathtt{o3 : MonomialIdeal of R}$
\medskip

\noindent $\mathtt{i5 : M = R^1/SRideal} ;$
\medskip

\noindent $\mathtt{o5 = cokernel | x1x2 x3x4 |}$
\medskip

\noindent $\mathtt{o5 : R-module, quotient of R}$
\medskip

\noindent $\mathtt{i6 : Mres = res\, \, M} ;$
\medskip

\noindent $\mathtt{o6 = R^1 \quad \longleftarrow \quad R^2 \quad \longleftarrow \quad R^1 \quad \longleftarrow \quad 0}$ 

\hspace{2mm} $\mathtt{ 0}$ \hspace{15mm} $\mathtt{1}$\hspace{16mm}$\mathtt{2}$ \hspace{14mm} $\mathtt{3}$
\medskip

\noindent $\mathtt{o6 : ChainComplex}$
\medskip

\noindent $\mathtt{i7 : \, peek\, \, \, betti \, \, \,  Mres }$
\medskip

\noindent $\mathtt{o7 = BettiTally{(0, \lbrace 0, 0, 0, 0, 0\rbrace , 0) => 1}}$ 

\hspace{21mm} $\mathtt{(1, \lbrace 1, 1, 0, 0, 0\rbrace , 2) => 1}$
 
\hspace{21mm} $\mathtt{(1, \lbrace 0, 0, 1, 1, 0\rbrace, 2) => 1}$
 
\hspace{21mm} $\mathtt{(2, \lbrace 1, 1, 1, 1, 0\rbrace , 4) => 1}$
\\

The last output labelled BettiTally has four columns. It starts with the level of the resolution. Then it follows the representation of $\sigma$ as an array of 0s and 1s between two braces. The third column indicates the cardinality of $\sigma$, and finally the last number indicates the multiplicity of $\sigma$.  In this example we can see that $\sigma= \lbrace 1, 1, 1, 1, 0 \rbrace$ at level 2 leads to $d_L(\Delta) = \max_{i, \sigma} \lbrace R_{i,\sigma} \rbrace = 2$.

%
% COMPLETE BIPARTITE GRAPH
%

\begin{example}{\bf(Complete Bipartite Graph)}
This example shows $d_L = \dhom > d_H$.  Let $\Delta(\C) = K_{r,r}$,  be the complete bipartite graph on $2r$ vertices.  Set $r=4$ as in Figure \ref{fig:bigraph}.  
The corresponding Stanley-Reisner ideal is 
$$I_{\Delta} = \langle x_1 x_2, x_1 x_3, x_1 x_4,
x_2 x_3, x_2 x_4, 
x_3 x_4, x_5 x_6, x_5 x_7, x_5 x_8,
x_6 x_7, x_6 x_8,  x_7 x_8,   \rangle$$

\begin{figure}[!h]
\begin{center}
\includegraphics[width=2.5in]{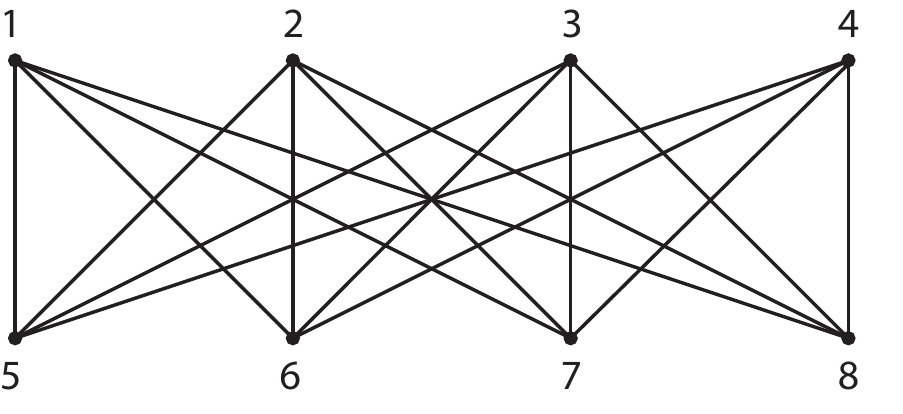}
\end{center}
\caption{Complete bipartite graph on $2r$ vertices for $r=4$.  The Helly dimension $d_H(\Delta)=1$, while the Leray dimension is $d_L(\Delta)=2$.}
\label{fig:bigraph}
\end{figure}

The minimal free resolution using the Stanley-Reisner ideal can be computed with {\tt Macaulay2}. We obtain the following resolution
$$0 \longleftarrow S\slash I_{\Delta}  \longleftarrow F_0  \longleftarrow F_1^{12}  \longleftarrow F_2^{52} \longleftarrow F_{3}^{102} \longleftarrow F_4^{100} \longleftarrow F_5^{48} \longleftarrow F_6^{9}  \longleftarrow0   $$
The calculation shows that at the first level of the resolution we have twelve elements $\sigma_1, \dots \sigma_{12}$ in $F_{1}^{12}$ with length $|\sigma_k| = 2$ for $k=1, \dots, 12$.  Thus, $d_H(\Delta) = \max_{\sigma} \lbrace R_{1,\sigma} \rbrace = 1$, as expected since $\Delta$ is a clique complex. For $1 < i \leq 5$, we have $\max_{i, \sigma} \lbrace R_{i,\sigma} \rbrace =2$.  Finally, at step $i = 6$, we find one $\sigma = 11111111$ with $|\sigma| = 8$ giving  $R_{i, 1\cdots 1}  = 2$.  Consequently, it follows from Theorem \ref{thm:main} that
\begin{itemize}
\item[]  $ d_L(\Delta) = \max_{i, \sigma} \lbrace R_{i,\sigma} \rbrace = 2,$

\item[] $d_H(\Delta) = \max_{\sigma}  \lbrace R_{1, \sigma} \rbrace = 1,$

\item[] $\dhom(\Delta) = \max_{i} \lbrace R_{i, 1\cdots 1} \rbrace = 2.$ 
\end{itemize}

%Thus, $R_{1, \sigma} = 1$.  At the higher levels $F_i$ for $5\geq i \geq1$ we find that the length of the elements is of two types. At each $F_i$ for $5 \geq i \geq 1$, we find elements $\sigma \in F_i$ whose difference between its length $|\sigma|$ and the step $i$ of the resolution is either 1 or 2. Denote by $\tilde{\sigma}$ those elements such that $|\tilde{\sigma}| - i  = 2$. The Betti numbers are 
%$$ \beta_{i, \tilde{\sigma}} = 1 \quad \textnormal{ with } \quad i=2, |\tilde{\sigma}| = 4, \quad i=3, |\tilde{\sigma}| = 5, \quad i = 4, |\tilde{\sigma}| = 6, \quad i = 5, |\tilde{\sigma}| = 7, \quad i= 6, |\tilde{\sigma}| = 8,$$
%and this leads to $R_{i, \tilde{\sigma}} = 2$ for $1 \leq i \leq 5$. 
%
% 
\end{example}
%

%
%%%%%%%%%%%%%%%%%%%%%%%%%%%%%%%%%%%%%%%%%%%
%
%               BIBLIOGRAPHY
%
%%%%%%%%%%%%%%%%%%%%%%%%%%%%%%%%%%%%%%%%%%%
%

\end{document}